\documentclass[11pt,twoside,a4paper,preprint,pre,floats,aps,amsmath,amssymb]{article}       

\usepackage[utf8]{inputenc} 
\usepackage[english]{babel}

\usepackage{amsthm}

\usepackage[left=1.5in,right=1.3in,top=1in,bottom=1.2in,includefoot,includehead,headheight=13.6pt]{geometry} 

\newcommand{\D}{\displaystyle}            

\usepackage{graphicx} 

\RequirePackage{fix-cm}

\usepackage{booktabs} 
\usepackage{array} 
\usepackage{paralist} 
\usepackage{verbatim} 
\usepackage{subfig} 
\newtheorem{theorem}{Theorem}[section]
\newtheorem{corollary}{Corollary}[section]
\newtheorem{lemma}{Lemma}[section]

\newtheorem{ex}{Exemple}[section]
\newtheorem{defn}{Definition}[section]
\newtheorem{rem}{Remarks}[section]
\newtheorem{prop}{Proposition}[section]

\newtheorem*{dem}{Proof}
\usepackage{amsthm}
\usepackage{amssymb}

\usepackage{fancyhdr} 
\usepackage{sectsty}
\allsectionsfont{\sffamily\mdseries\upshape} 

\usepackage[nottoc,notlof,notlot]{tocbibind} 
\usepackage[titles,subfigure]{tocloft} 

\usepackage{minitoc}
\usepackage{aecompl}

\usepackage[intoc]{nomencl}

\makenomenclature

\graphicspath{{.}{images/}}

\usepackage{color}
\definecolor{linkcol}{rgb}{0.87,0.00,0.87}
\definecolor{citecol}{rgb}{0.50,0.00,1.00}
\usepackage{fancyhdr}
\usepackage{colortbl}
\arrayrulecolor{black}

\usepackage{hyperref}
\hypersetup{
pdftex, 											
bookmarks = true, 		
bookmarksopen=true,         
bookmarksnumbered=true,    
breaklinks=true, 
pdfpagelayout=SinglePage,
pdfpagemode=UseOutlines,   
pdfstartview=FitH,              
colorlinks = true,	 	
urlcolor = blue, 		
pdfborder = {0 0 0}, 	
linkcolor=magenta,              
urlcolor=magenta,                     
anchorcolor=black,              
citecolor=blue,                 
frenchlinks=true,                
pdffitwindow=true,             
pdftoolbar=true,                 
pdfmenubar=true,               
pdfauthor={},     	             
pdftitle={},    	             
pdfkeywords={gggggggggg}, 						
pdfcreator=PDFLaTeX,         %
pdfproducer=PDFLaTeX,       %
}

\providecommand{\keywords}[1]{\textbf{\textit{Keywords:}} #1}
\usepackage[affil-it]{authblk}

\usepackage{amsmath}
\numberwithin{equation}{section}

\title{$q$-Difference Equations Associated with the Rubin’s $q$-Difference Operator $\partial_q$}
\newcommand*{\email}[1]{{E-mail: #1}}
\author[,1]{Meniar Haddad\thanks{\email{\texttt{meniar.haddad@fst.rnu.tn}}}}
\author[,2]{Marwa Mastouri\thanks{\email{\texttt{marwamastouri93@gmail.com}}}}

\affil[1]{Department of Mathematics, Faculty of Sciences of Bizerte, 7021 Zarzouna, Tunisia}
\affil[2]{Department of Mathematics, Faculty of Sciences of Bizerte, 7021 Zarzouna, Tunisia}

\begin{document}
\maketitle

\begin{abstract}
The aim of this paper is to prove the existence and uniqueness of solutions of the following $q$-Cauchy problem of second order linear $q$-difference problem associated with the Rubin's $q$-difference operator $\partial_q$ in a neighborhood of zero
\begin{equation}
\left\{
  \begin{array}{cc}
      q\,a_0(x)\,\partial_q^2y(qx)\,+\,a_1(x)\,\partial_qy(x)\,+\,a_2(x)y(x)&\;=\;b(x),\quad \hbox{if $y$ is odd;}\\
      q\,a_0(x)\partial_q^2y(qx)\,+\,q\,a_1(x)\partial_qy(qx)\,+\,a_2(x)y(x)&\;=\;b(x),\quad \hbox{if $y$ is even,}
 \end{array}
\right.
    \end{equation}
with the initial conditions
\begin{equation}
   \partial_{q}^{i-1}y(0)= b_{i};\quad b_{i} \in{\mathbb{C}},\; i=1,2
 \end{equation}
where $a_i$, $i=0,1,2$, and $b$ are defined, continuous at zero and bounded on an interval $I$ containing zero such that $a_0(x)\neq 0$ for all $x\in I$.
Then, as application of the main results, we study the second order homogenous linear $q$-difference equations as well as the $q$-Wronskian associated with the Rubin's $q$-difference operator $\partial_q$. Finally, we construct a fundamental set of solutions for the second order linear homogeneous $q$-difference equations in the cases when the coefficients are constants and $a_1(x)=0$ for all $x\in I$.\\

\end{abstract}
\section{Introduction}

Studies on $q$-difference equations appeared already at the beginning of the last century in intensive works especially by F.H. Jackson$~{\cite{N}}$, R.D. Carmichael$~{\cite{R.D}}$, T.E. Mason$~{\cite{T.E}}$, C.R. Adams$~{\cite{C.R}}$ and other authors such us Poincare, Picard, Ramanujan. Apart from this old history of $q$-difference equations, the subject received a considerable interest of many mathematicians and from many aspects, theoretical and practical.

Since years eighties$~{\cite{Hahn}}$, an intensive and somewhat surprising interest in the subject reappeared in many areas of mathematics and applications including mainly new $q$-difference calculus and $q$-orthogonal polynomials, $q$-combinatorics, $q$-arithmetics, $q$-integrable systems.

The present article is devoted for developing the theory of the second order linear $q$-difference equations associated with the Rubin's $q$-difference operator $\partial_q$ in a neighborhood of zero, studying the $q$-wronskian associated with the Rubin’s $q$-difference operator $\partial_q$ and to showing how it plays a central role in solving the second order linear $q$-difference equations. As application, we study the second order linear homogeneous $q$-difference equations associated with the Rubin’s $q$-difference operator $\partial_q$ when the coefficients are constants. As M.H. Annaby and Z.S. Mansour$~{\cite{j6}}$, we establish some results associated with the Rubin’s $q$-difference operator $\partial_q$.

We mention that, in this paper, we will follow all the so-mentioned works which, for a convergence argument, they imposed that the parameter $q\in]0,1[$ satisfies the condition
$$\frac{\ln{(1-q)}}{\ln(q)}\in 2\,\mathbb{Z}.$$

This paper is organized as follows: in Section 2, we recall some necessary fundamental concepts of quantum analysis. Section 3 is devoted to prove the existence and uniqueness of solutions of the $q$-Cauchy problem of second order linear $q$-difference equations in a neighborhood of zero. In Section 4, we study the second order linear $q$-difference equations associated with the Rubin’s $q$-difference operator $\partial_q$. In Section 5, we introduce the $q$-Wronskian associated with the Rubin’s $q$-difference operator $\partial_q$ and we establish some of its properties. In Section 6, we construct a fundamental set of solutions for the second order linear homogeneous $q$-difference equations when the coefficients are constants. Finally, in Section 7, we purpose some examples.
\section{Notations and preliminaries}
In this section, we introduce some necessary fundamental concepts of quantum analysis which will be used in this paper. For this purpose, we refer the reader to the book by G. Gasper and M. Rahman $~{\cite{L}}$, for the definitions, notations and properties of the $q$-shifted factorials and the $q$-hypergeometric functions.
Throughout this paper, we assume $q \in]0, 1[$ and we denote
$${\mathbb R}_{q}\; =\; \{\pm q^k , k\in \mathbb Z\},\quad {\mathbb R}_{q,+} \;=\; \{+ q^k , k\in \mathbb Z\}\quad \text{and}\quad \widetilde{{\mathbb R}}_{q}\;=\; {\mathbb R}_{q}\, \bigcup\, \{0\}.$$
\subsection{Basic symbols}
For a complex number $a$, the $q$-shifted factorials are defined by:
$${(a,q)}_{0}\;:=\;1;\quad{(a,q)}_{n}\;:=\;\displaystyle\prod_{k=0}^{n-1}(1-a{q}^{k}),\; n\geq 1\quad \text{and}\quad {(a,q)}_{\infty}\;=\;\prod_{k=0}^{\infty}(1-a{q}^{k}).$$
We also denote
$$ {[n]}_{q}\;=\;\frac{1-q^n}{1-q}\quad\text{and}\quad{[n]_{q}}!\;=\;\frac{(q,q)_{n}}{(1-q)^n},\quad n\in\mathbb{N}.$$
If we change $q$ by $q^{-1}$, we obtain
$$[n]_{q^{-1}}!\;=\;q^{-\frac{n(n-1)}{2}}\, [n]_{q},\quad n\in\mathbb{N}.$$
The Gauss $q$-binomial coefficient $\displaystyle \left[{\begin{matrix}n\\k\end{matrix}}\right]_{q}$ is defined by (see$~{\cite{L}}$)
$${\displaystyle \left[{\begin{matrix}n\\k\end{matrix}}\right]_{q}}\;=\;\D{{\frac {[n]_{q}!}{[n-k]_{q}![k]_{q}!}}} \;=\;\frac{(q,q)_{n}}{(q,q)_{k}\,(q,q)_{n-k}},\quad n\geq k\geq 0.$$
Moreover, the $q$-binomial theorem is given by
$$(-a,q)_n\;=\;\D\sum_{k=0}^{n}{\displaystyle \left[{\begin{matrix}n\\k\end{matrix}}\right]_{q}}q^{\frac{k(k-1)}{2}}a^k,\quad a\in\mathbb{C}.$$
\subsection{Operators and elementary $q$-special functions}
Let $\mathcal{A} \subseteq \mathbb{R}$ be a $q$-geometric set containing zero and satisfying for every $x\in\mathcal{A}$, $\pm q^{\pm1}x\in\mathcal{A}$. A $q$-difference equation is an equation that contains $q$-derivatives of a function defined on $\mathcal{A}$. Let $f$ be a function, real or complex valued, defined on a $q$-geometric set $\mathcal{A}$. The $q$-difference operator $D_q$, the Jackson $q$-derivative is defined by
\begin{equation}\label{ddd}
D_{q}f(x)\;:=\;\frac{f(x)-f(qx)}{(1-q)x},\quad \text{for all}\;x\in\mathcal{A}\backslash\{0\}.
  \end{equation}
In the $q$-derivative, as $q\rightarrow 1$, the $q$-derivative is reduced to the classical derivative.\\
The $q$-derivative at zero is defined by
\begin{equation}
D_qf(0)\;:=\;\lim_{n\rightarrow +\infty}\dfrac{f(q^nx)-f(0)}{q^nx},\quad \text{for all}\;x\in\mathcal{A},
\end{equation}
if the limit exists and does not depend on $x$.\\
We recall that for a function $f(x)$ defined on $\mathcal{A}$ and $n\in\mathbb{N}$, we have
\begin{equation}\label{lm}
\displaystyle D_{q}^{n}f(x)\;=\;{\frac {(-1)^{n}q^{-\frac{n(n-1)}{2}}}{(1-q)^{n}x^{n}}}\,\sum _{k=0}^{n}(-1)^{k}\,{\displaystyle \left[{\begin{matrix}n\\k\end{matrix}}\right]_{q}}\,q^{\frac{k(k-1)}{2}}f(q^{n-k}x).
\end{equation}
The Rubin's $q$-difference operator $\partial_q$ is defined in $~{\cite{&,3+}}$ by
\begin{equation}\label{/}
\partial_{q}f(x)\;=\;\left\{
                         \begin{array}{ll}
                           \D\dfrac{f(q^{-1}x)+f(-q^{-1}x)-f(qx)+f(-qx)-2f(-x)}{2(1-q)x}&,\;  \hbox{$x\neq0$;} \\
                           \D\lim_{x\rightarrow0}\partial_{q}f(x)&,\;  \hbox{$x=0$.}
                         \end{array}
                       \right.
\end{equation}
It is straightforward to prove that if a function $f$ is differentiable at a point $z$, then
$$\D{\lim_{q\rightarrow 1^{-}}\partial_{q}f(x)\,=\, f\,'(x)}.$$
A repeated application of the Rubin's $q$-difference operator $n$ times is denoted by:
$$\partial_{q}^0f\;=\;f\quad{\text{and}}\quad \partial_{q}^{n+1}f\;=\;\partial_{q}(\partial_{q}^{n}f),\quad n\in\mathbb{N}.$$
We define the $q$-shift operators by:
\begin{eqnarray}
(\Lambda_qf)(x) \;=\,& f(qx)\quad &\text{and}\quad (\Lambda_q^{-1} f)(x)\;=\; f(xq^{-1}),\label{Lam}\\
\partial_q\,\Lambda_qf\;=\,& q\Lambda_q\partial_q \quad &\text{and}\quad
\partial_q\,\Lambda_q^{-1}f\;=\; q^{-1}\,\Lambda_q^{-1}\,\partial_q.\label{Lam1}
\end{eqnarray}
We mention that $\partial_q$ is closely related to the Jackson’s $q$-derivative operator $D_q$ and use relation $~(\ref{lm})$, we can easily prove the following result:
\begin{prop}{$:$}\label{code}
Let $f$ be a function defined on $\mathbb{R}_q$. Then for all $n\in\mathbb{N}$ we have:
\begin{enumerate}
\item \quad \begin{equation}\label{seconde}
\partial_{q}^{2n}f\;=\;q^{-n(n+1)}(D_{q}^{2n}f_{e})o\Lambda_{q}^{-n}+q^{-n^2}(D_{q}^{2n}f_{o})o\Lambda_{q}^{-n},\end{equation}
\item \quad \begin{equation}
\partial_{q}^{2n+1}f\;=\;q^{-(n+1)^2}(D_{q}^{2n+1}f_{e})o\Lambda_{q}^{-(n+1)}+q^{-n(n+1)}(D_{q}^{2n+1}f_{o})o\Lambda_{q}^{-n},\end{equation}
\end{enumerate}
where $f_{e}$ and $f_{o}$ are, respectively, the even and the odd parts of $f$ and $\Lambda_{q}^{-n}$ is the function given by $\Lambda_q^{-n}(x)\,=\,q^{-n}x$.
\end{prop}

 A right inverse to the $q$-derivative, the $q$-integration is given by Jackson  as
 \begin{equation}
\int_{0}^{x}f(t)d_{q}t \; := \; x\,(1-q)\,\displaystyle\sum_{n=0}^{+\infty}q^{n}f(q^{n}\,x),\quad \text{for all}\;x\in\mathcal{A},
\end{equation}
provided that the series converges. In general,
\begin{equation}
  \D\int_{a}^{b}f(t)d_{q}t \; :=\; \int_{0}^{b}f(t)d_{q}t\, -\, \int_{0}^{a}f(t)d_{q}t, \quad \text{for all}\;a,b\in\mathcal{A}.
\end{equation}
The $q$-integration for a function is defined in $~{\cite{HW}}$ by the formulas
\begin{align}
  \D \int_{0}^{+\infty}f(t)d_{q}t &\;=\; (1-q)\,\displaystyle\sum_{n=-\infty}^{+\infty}\,q^{n}\,f(q^{n}),\label{a29} \\
    \D \int_{-\infty}^{0}f(t)d_{q}t &\;=\; (1-q)\,\displaystyle\sum_{n=-\infty}^{+\infty}\,q^{n}\,f(-q^{n}),\label{a29*} \\
\D\int_{-\infty}^{+\infty}f(t)d_{q}t &\;=\; (1-q)\,\D\sum_{n=-\infty}^{+\infty}\,q^{n}\,\left[f(q^{n})\,+\,f(-q^{n})\right],
\end{align}
provided that the series converges.
\begin{rem}{$:$}
Note that when $f$ is continuous on $[0,a]$, it can be shown that
$$\D\lim_{q\rightarrow1}\int_{0}^{a}f(t)d_{q}t \;=\; \int_{0}^{a}f(t)dt.$$
\end{rem}
\begin{defn}
A function $f$ which is defined on $\mathcal{A}$, $0\in\mathcal{A}$, is said to be $q$-regular at zero if
\begin{equation}
  \lim_{n\rightarrow \infty}f(xq^n)\;=\;f(0),\quad \text{for all}\;x\in\mathcal{A}.
\end{equation}
\end{defn}
Through the remainder of the paper, we deal only with functions that are $q$-regular at zero.
The following results hold by direct computation.
\begin{lemma}\label{m"}
  Let $f$ be a function defined on $\mathcal{A}$. Then, if $\D\int_{0}^{x}f(t)\,d_{q}t$ exists, we have:
  \begin{enumerate}
\item \quad for all integer $n$,\quad $\D\int_{-\infty}^{+\infty}f(q^nt)d_{q}t\;=\;q^{-n}\,\D\int_{-\infty}^{+\infty}f(t)d_{q}t.$
\item \quad if $f$ is odd, then
\begin{equation}\label{mr'}
\partial_{q}\D\int_{0}^{x}\,f(t)d_{q}t  \;=\; f(x),
\end{equation}
and
\begin{equation}\label{mr}
\D\int_{0}^{x}\,\partial_{q}f(t)d_{q}t  \;=\; f(x)\,-\, \lim_{n\rightarrow \infty}f(xq^n).
\end{equation}
\item \quad if $f$ is even, then
\begin{equation}\label{mr1'}
\partial_{q}\D\int_{0}^{x}\,f(t)d_{q}t  \;=\; f(q^{-1}x),
\end{equation}
and
\begin{equation}\label{mr1}
\D\int_{0}^{x}\,\partial_{q}f(t)d_{q}t  \;=\; f(q^{-1}x)\,-\,\lim_{n\rightarrow \infty}f(xq^{n-1}).
\end{equation}
\end{enumerate}
\end{lemma}
\begin{prop}{\bf[The rule of $q$-integration by parts]}
Let $f$ and $g$ be two functions defined on $[-a,a]$, for all $a>0$. Then, if $\D\int_{-a}^{a}\partial_qf(t)\,g(t)\,d_{q}t$ exists, the rule of $q$-integration by parts is given by:
\begin{equation}
\D\int_{-a}^{a}\partial_qf(t)\,g(t)\,d_{q}t\;=\;2\left[f_e(q^{-1}a)g_o(a)\,+\,f_o(a)g_e(q^{-1}a)\right]\,-\,\D\int_{-a}^{a}f(t)\,\partial_qg(t)\,d_{q}t,
\end{equation}
where $f_e$, $f_o$ are the even and the odd parts of $f$ respectively and $g_e$, $g_o$ are the even and the odd parts of $g$ respectively.
\end{prop}
\begin{corollary}
  Let $f$ and $g$ be two functions defined on $\mathbb{R}_q$. Then, if $\D\int_{-\infty}^{+\infty}\partial_qf(t)\,g(t)\,d_{q}t$ exists, the rule of $q$-integration by parts is given by:
\begin{equation}
  \D\int_{-\infty}^{+\infty}\partial_qf(t)\,g(t)\,d_{q}t\;=\;-\,\D\int_{-\infty}^{+\infty}f(t)\,\partial_qg(t)\,d_{q}t.
\end{equation}
\end{corollary}
\begin{defn}
Let $f$ be a function defined on $\mathcal{A}$. We say that $f$ is $q$-integrable on $\mathcal{A}$ if and only if $\D \int_{0}^{x}f(t)d_{q}t$ exists for all $x\in\mathcal{A}$.
\end{defn}

The $q$-trigonometric functions $q$-cosine and $q$-sine (see$~{\cite{&,3+}}$) are defined, respectively, on $\mathbb{C}$ by:
\begin{equation}\label{cos}
\cos(x,q^2)\;:=\;\D\sum_{n=0}^{+\infty}\,(-1)^n\,b_{2n}(x,q^2)\end{equation}
and
\begin{equation}\label{sin}
\sin(x,q^2)\;:=\;\D\sum_{n=0}^{+\infty}\,(-1)^n\,b_{2n+1}(x,q^2)
\end{equation}
where
\begin{equation}\label{mas}
b_n(x,q^2)\;=\;q^{[\frac{n}{2}]([\frac{n}{2}]+1)}\dfrac{x^{n}}{[n]_q!},\quad n\in\mathbb{N}.
\end{equation}
These two functions induce a $\partial_q$-adapted $q$-analogue exponential function as
\begin{equation}\label{exp3}
e(x,q^2)\;:=\;\cos(-ix,q^2)\,+\,i\sin(-ix,q^2)\;=\;\D\sum_{n=0}^{+\infty}b_{n}(x,q^2).
\end{equation}
Remark that $e(x,q^2)$ is absolutely convergent for all $x$ in the complex plane since both of its component functions are. Moreover, $\D\lim_{q\rightarrow1^{-}}e(x,q^2)\,=\,e^x$ (exponential function) point-wise and uniformly on compacts. The following results hold by direct computation.
\begin{lemma}\label{sit}
\begin{enumerate}
  \item  For all $x\in\mathbb{C}$ and $\lambda\in\mathbb{C}$, we have
\begin{equation}
\partial_q\cos(\lambda x,q^2)\;=\;-\,\lambda\,\sin(\lambda x,q^2),\quad \partial_q\sin(\lambda x,q^2)\;=\;\lambda\,\cos(\lambda x,q^2)
\end{equation}
and
\begin{equation}
  \partial_q e(\lambda x,q^2)\;=\;\lambda\,e(\lambda x,q^2).
\end{equation}
  \item  For all function $f$ defined on $\mathcal{A}$, we have
\begin{equation}\label{PIM}
\partial_{q}f(x)\;=\;\D\dfrac{f_e(q^{-1}x)-f_e(x)}{(1-q)x}\,+\,\D\dfrac{f_o(x)-f_o(qx)}{(1-q)x},\quad x\in\mathcal{A}\backslash\{0\}\end{equation}
Here, $f_e$ and $f_o$ are the even and the odd parts of $f$, respectively.
  \item  For two functions $f$ and $g$ defined on $\mathcal{A}$, we have  \\
  -- \quad if $f$ is even and $g$ is odd then
\begin{equation}\label{marwa1}
  \partial_{q}(fg)(x)\;=\;f(x)\partial_{q}g(x)+q\,g(qx)\partial_{q}f(qx)\;=\;q\,g(x)\partial_{q}f(qx)+f(qx)\partial_{q}g(x).
\end{equation}
  -- \quad if $f$ and $g$ are even then
\begin{equation}\label{marwa2}
  \partial_{q}(fg)(x)\;=\;g(q^{-1}x)\partial_{q}f(x)+f(x)\partial_{q}g(x).
\end{equation}
-- \quad if $f$ and $g$ are odd then
\begin{equation}\label{marwa3}
  \partial_{q}(fg)(x)\;=\;q^{-1}\,g(q^{-1}x)\partial_{q}f(q^{-1}x)+q^{-1}f(x)\partial_{q}g(q^{-1}x).
\end{equation}
\end{enumerate}
\end{lemma}
\section{$q$-Initial Value Problems in a Neighborhood of Zero}
In this section, we prove the existence and uniqueness of solutions of the $q$-Cauchy problem of second order $q$-difference equations in a neighborhood of zero. 
In the sequel $X$ is a Banach space with norm $\|\cdot\|$ and $I\subseteq\mathbb{R}$ is an interval containing zero.
\begin{defn}
Let $S$ and $R$ be defined, respectively, by
 $$S(y_0,\beta)\,:=\,\{y\in\mathbb{X}\,:\,\|y-y_0\|\,\leq\,\beta\},$$
 and
 $$R\,:=\,\{(x,y)\in I \times\mathbb{X}\,:\,|x|\leq \alpha,\,\|y-y_0\|\,\leq\,\beta\},$$
 where $y_0\in\mathbb{X}$ and $\alpha$, $\beta$ are fixed positive real numbers.
\end{defn}
By a $q$-initial value problem ($q$-IVP) in a neighborhood of zero we mean the problem of finding continuous functions at zero 
satisfying system \begin{equation}\label{IP1}
\partial_{q}{y(x)}\;=\;f(x,y(x)),\quad y(0)\;=\;y_0,\quad {x}\in I.\end{equation}
\begin{theorem}
Let $f:R\rightarrow {\mathbb{X}}$ be a continuous function at $(0,y_0)\in R$ , and $\phi$ be a function defined on $I$. Then $\phi$ is a solution of the
$q$-IVP $~(\ref{IP1})$ if, and only if,
\begin{enumerate}
  \item For all $x\in I$, $(x,\phi(x))\in R$.
  \item $\phi$ is continuous at zero.
  \item For all $x\in I$,
  \begin{equation}\label{IP5}
\int_{0}^{x}f(t,\phi(t))d_qt\;=\;\left\{
  \begin{array}{cc}
\phi(x)\,-\,y_0,&\quad \hbox{if $y$ is odd;}\\
\phi(q^{-1}x)\,-\,y_0,&\quad \hbox{if $y$ is even,}
\end{array}
\right.
\end{equation}
\end{enumerate}
\end{theorem}
\begin{dem}
Let $\phi$ be a solution of the $q$-IVP $~(\ref{IP1})$. Then
\begin{equation}\label{IP2}
\partial_{q}{\phi(x)}\;=\;f(x,\phi(x)),\quad \phi(0)\;=\;y_0,\quad {x}\in I.\end{equation}
which implies $(x,\phi(x))\in R$ for all $x \in I$. Also, since $\phi$ is $q$-differentiable on $I$, then it is continuous at zero. Finally, integrating both sides of $~(\ref{IP2})$ from 0 to $x$ and using Lemma$~\ref{m"}$ we get
$$\int_{0}^{x}f(t,\phi(t))d_qt\;=\;\left\{
  \begin{array}{cc}
\phi(x)\,-\,\phi(0),&\quad \hbox{if $y$ is odd;}\\
\phi(q^{-1}x)\,-\,\phi(0),&\quad \hbox{if $y$ is even,}
\end{array}
\right.$$
Consequently,
$$\int_{0}^{x}f(t,\phi(t))d_qt\;=\;\left\{
  \begin{array}{cc}
\phi(x)\,-\,y_0,&\quad \hbox{if $y$ is odd;}\\
\phi(q^{-1}x)\,-\,y_0,&\quad \hbox{if $y$ is even,}
\end{array}
\right.$$
Conversely, assume the items (1), (2) and (3) are satisfied, then $\phi$ is ordinary differentiable at zero. Consequently, it is $q$-differentiable on $I$ with $\partial_{q}{\phi(x)}\;=\;f(x,\phi(x))$ and $\phi(0)\;=\;y_0$. Therefore, $\phi$ is a solution of the $q$-IVP $~(\ref{IP1})$.
\hfill$\blacksquare$
\end{dem}
To prove the existence and uniqueness of the solution of the $q$-IVP $~(\ref{IP1})$, we need some preliminary results:
\begin{defn}
 Let ${\displaystyle (\mathbb{X},d)}$ be a complete metric space. Then a map ${\displaystyle T\colon \mathbb{X}\to \mathbb{X}}$ is called a contraction mapping on ${\displaystyle \mathbb{X}}$ if there exists $0\leq\rho<1$ such that
$${\displaystyle d(T(x),T(y))\leq \rho\,d(x,y)}$$
for all ${\displaystyle x,y}\in {\displaystyle \mathbb{X}}$.
\end{defn}
\begin{theorem}{\bf(Banach's fixed point theorem).}\label{BPT}\;Let ${\displaystyle (\mathbb{X},d)}$ be a non-empty complete metric space with a contraction mapping ${\displaystyle T\colon \mathbb{X}\to \mathbb{X}}$. Then $T$ has a unique fixed-point $x^* \in \mathbb{X}$ (i.e. $T(x^*) = x^*$). Furthermore, for any $x\in \mathbb{X}$ and $n\geq1$ the iterative sequence $\{T^n(x)\}$ converges to $x^*$.
\end{theorem}
\begin{theorem}\label{th1}
Let $f:R\rightarrow {\mathbb{X}}$ be a continuous function at $(0,y_0)\in R$ and satisfies the Lipschtiz condition with respect to $y$ in $R$, that is, there exists a positive constant $L$ such that:
\begin{equation}
\big\| f(x,y_{1})-f(x,y_{2})\big\| \leq{L}\big\|y_{1}-y_{2} \big\|, \quad\textit{for all } (x,y_{1}), (x,y_{2})\in{R}.
\end{equation}
Then the $q$-IVP $~(\ref{IP1})$ has a unique solution on $[-h,h]$, where
$$h\;:=\;\min\{\alpha,\frac{\beta}{L\beta\,+\,M},\frac{\rho}{L}\}\quad \text{with}\quad M\,=\,\sup_{(x,y)\in{R}}\|f(x,y)\|<+\infty,\quad 0\leq\rho<1.$$
\end{theorem}
\begin{dem}
We prove the theorem for $x\in[0,h]$ and the proof for $x \in [-h,0]$ is similar. Define the operator $T$ by
\begin{equation}\label{IP3}
  T\,y(x)\;=\;y_0\,+\,\int_{0}^{x}f(t,y(t))d_qt.
\end{equation}
Let $\mathcal{C}([0,h])$ be the space of all continuous functions at zero and bounded on the interval $[0,h]$ with the supremum norm such that for $y\in\mathcal{C}([0,h])$, we have
$$\|y\|_{\infty}\;=\;\displaystyle\sup_{x\in[0,h] }\,\|y(x)\|.$$
 This space is complete. Let
 $$\widetilde{S}(y_0,\beta)\,:=\,\{y\in\mathcal{C}([0,h])\,:\,\|y\,-\,y_0\|_{\infty}\,\leq\,\beta\},\quad \beta >0.$$
As $\widetilde{S}(y_0,\beta)$ is a closed subset of the complete space $\mathcal{C}([0,h])$, it is also complete metric space. 
 First, we prove that $T:\, \widetilde{S}(y_0,\beta)\rightarrow \widetilde{S}(y_0,\beta)$. Let $\phi\in \widetilde{S}(y_0,\beta)$,
 \begin{align*}
 \big\|T\,\phi(x)\,-\,y_0\big\|_\infty&\;=\;\displaystyle\sup_{x\in[0,h] }\,\big\| \int_{0}^{x}f(t,\phi(t))d_qt\big\|\\
 &\;=\;\displaystyle\sup_{x\in[0,h] }\,\big\| \int_{0}^{x}\left(f(t,\phi(t))\,-\,f(t,y_0)\,+\,f(t,y_0)\right)d_qt\big\|\\
 &\;\leq\;\displaystyle\sup_{x\in[0,h] }\, \int_{0}^{x}\big\|f(t,\phi(t))\,-\,f(t,y_0)\,+\,f(t,y_0)\big\|d_qt\\
 &\;\leq\; \displaystyle\sup_{x\in[0,h] }\,\int_{0}^{x}\big\|f(t,\phi(t))\,-\,f(t,y_0)\big\|\,+\,\big\|f(t,y_0)\big\|d_qt\\
  &\;\leq\;\displaystyle\sup_{x\in[0,h] }\, \int_{0}^{x}\left(L\,\big\|\phi(t)\,-\,y_0\big\|\,+\,M\right)d_qt\\
  &\;\leq\;\left(L\,\big\|\phi\,-\,y_0\big\|_\infty\,+\,M\right)\,\int_{0}^{x}d_qt\\
 &\;\leq\; \left(L\,\beta\,+\,M\right)\,\int_{0}^{x}d_qt\\
 &\;\leq\; \left(L\,\beta\,+\,M\right)\,x\\
 &\;\leq\; \left(L\,\beta\,+\,M\right)\,h.\end{align*}
 Then, using the fact that $h\leq \frac{\beta}{L\beta\,+\,M}$, we obtain
 $$\|T\,\phi\,-\,y_0\|_\infty\leq \beta,\quad \text{i.e.},\quad T\,\phi\in \widetilde{S}(y_0,\beta).$$
  Next, we prove that $T$ is a contraction mapping. Assume that $\phi_1,\phi_2\in \widetilde{S}(y_0,\beta)$, then
\begin{align*}
 \big\|T\,\phi_1(x)\,-\,T\,\phi_2(x) \big\|&\;=\;\big\| \int_{0}^{x}\left(f(t,\phi_1(t))\,-\,f(t,\phi_2(t))\right)d_qt\big\|\\
 &\;\leq\; \int_{0}^{x}\big\|f(t,\phi_1(t))\,-\,f(t,\phi_2(t))\big\|d_qt\\
  &\;\leq\; \int_{0}^{x}L\,\big\|\phi_1(t)\,-\,\phi_2(t)\big\|d_qt\\
  &\;\leq\;L\,\big\|\phi_1\,-\,\phi_2\big\|_{\infty} \,\int_{0}^{x}d_qt\\
  &\;\leq\; L\,\big\|\phi_1\,-\,\phi_2\big\|_{\infty}\,x\\
  &\;\leq\;L\,\big\|\phi_1\,-\,\phi_2\big\|_{\infty} \,h.
 \end{align*}
 But, the fact that
 $$h\leq \frac{\rho}{L}\quad \text{with}\quad 0\leq\rho<1$$
 gives
 $$ \|T\,\phi_1(x)\,-\,T\,\phi_2(x) \|\;\leq\;\rho\,\|\phi_1\,-\,\phi_2\|_{\infty}.$$
 Then $T$ is a contraction mapping. By Banach’s fixed point theorem$~\ref{BPT}$, $T$ has a unique fixed point in $\widetilde{S}(y_0,\beta)$ and then the $q$-IVP $~(\ref{IP1})$ has a unique solution in $\widetilde{S}(y_0,\beta)$.
\hfill$\blacksquare$
\end{dem}
\begin{theorem}\label{th2}
Let the functions $f_{i}(x,y_{1},y_{2})$, $i=1,2$, be defined on $I \times\prod_{i=1}^{2} S_{i}(b_i, \beta_i)$, such that the following conditions are satisfied:
\begin{description}
  \item[(i)] for $y_{i}\in S_{i}(b_i,\beta_i)$, $i=1,2$, $f_{i}(x,y_{1},y_{2})$ are continuous at zero,
  \item[(ii)] there is a positive constant $L$ such that, for $x\in I$, $y_{i}, \tilde{y}_{i}\in S_{i}(b_i,\beta_i)$, $i=1,2,$ the following Lipschitz condition is satisfied:
      $$\big\| f_{i}(x,y_{1},y_{2})-f_{i}(x, \tilde{y}_{1},\tilde{y}_{2})\big\| \leq L \sum _{i=1}^{2}\|y_{i}-\tilde{y}_{i} \|.$$
\end{description}
Then there exists a unique solution of the $q$-initial value problem,
\begin{equation}\label{IVP4}
\partial_{q}y_{i}(x)=f_{i}\bigl(x,y_{1}(x),y_{2}(x) \bigr),\quad y_{i}(0)=b_i\in\mathbb{X},\; i =1,2,\; x \in I.
\end{equation}
\end{theorem}
\begin{dem}
Let $y_0=(b_{1},b_2)^{T}$ and $\beta=(\beta_{1},\beta_{2})^{T}$, where $(\cdot ,\cdot)^{T}$ stands for vector transpose. The function $f:I\times \prod_{i=1}^{2}S_{i}(b_i,\beta_i)\rightarrow{\mathbb{X}}\times{\mathbb{X}}$ is defined by
$$f(x,y_{1},y_{2})\;=\; (f_{1}(x,y_{1},y_{2}),f_{2}(x, y_{1}, y_{2}) )^{T}.$$
It is easy to show that system $~(\ref{IVP4})$ is equivalent to the $q$-IVP $~(\ref{IP1})$.\\
Since each $f_i$ is continuous at zero, $f$ is continuous at zero. We claim that $f$ satisfies the Lipschitz condition. Indeed for $y=(y_1,y_2)$ and $\tilde{y}=(\tilde{y}_1,\tilde{y}_2)$ in $\prod_{i=1}^{2}S_{i}(b_i,\beta_i)$
\begin{align*}
\big\| f(x,y)-f(x,\tilde{y}) \big\| &\;=\; \big\| f(x,y_{1},y_{2})-f(x, \tilde{y}_{1}, \tilde{y}_{2}) \big\| \\ &\;=\;\sum_{i=1}^{2} \big\| f_{i}(x,y_{1},y_{2})-f_{i}(x, \tilde{y}_{1},\tilde {y}_{2}) \big\| \\
&\;\leq\; L\,\sum_{i=1}^{2}\|y_{i}- \tilde{y}_{i}\|\\
&\;=\; L\,\|y-\tilde{y}\|.\end{align*}
Applying Theorem$~\ref{th1}$, there exists $h>0$ such that $~(\ref{IP1})$ has a unique solution on $[-h,h]$. Hence, the $q$-IVP $~(\ref{IVP4})$ has a unique solution on $[-h,h]$.
\hfill$\blacksquare$
\end{dem}
\begin{corollary}\label{cor1}
 Let $f(x,y_{1},y_{2})$ be a function defined on $I\times\D\prod_{i=1}^{2} S_{i}(b_i,\beta_i)$ such that the following conditions are satisfied:
 \begin{description}
 \item[(i)]  for any values of $y_{i}\in S_{i}(b_i,\beta_i)$, $i=1,2$, $f$ is continuous at zero,
\item[(ii)]  $f$ satisfies the following Lipschitz condition
$$\big\| f(x,y_{1},y_{2})-f(x,\tilde{y}_{1}, \tilde{y}_{2}) \big\| \leq L\,\sum_{i=1}^{2} \big\|y_{i} -\tilde{y}_{i}\big\|,$$\end{description}
where $L>0$, $y_{i},\tilde{y}_{i}\in S_{i}(b_i,\beta_i)$, $i=1,2$ and $x \in I$. Then
\begin{equation}\label{IP5}
\partial_{q}^{2}y(qx)\;=\;\left\{
  \begin{array}{cc}
f\bigl(x,y(x),\partial_{q}y(x)\bigr),\quad \hbox{if $y$ is odd;}\\
f\bigl(x,y(x),\partial_{q}y(qx)\bigr),\quad \hbox{if $y$ is even,}
\end{array}
\right.
\end{equation}
with the initial conditions
\begin{equation}
\partial_{q}^{i-1}y(0)=b_i,\quad i=1,2
\end{equation}
has a unique solution on $[-h,h]$.
\end{corollary}
\begin{dem}
Consider equation $~(\ref{IP5})$. It is equivalent to $~(\ref{IVP4})$, where $\{\phi_{i}(x)\}_{i=1}^{2}$ is a solution of $~(\ref{IVP4})$ if and only if $\phi_{1}(x)$ is a solution of $~(\ref{IP5})$. Here,
\begin{align*} f_{i}(x,y_{1},y_{2})\;=\;\left \{ \textstyle\begin{array}{l@{\quad}l} y_{2},& i=1, \\ f (x,y_{1},y_{2}),& i=2. \end{array}\displaystyle \right . \end{align*}
Hence, by Theorem$~\ref{th2}$, there exists $h>0$ such that system $~(\ref{IVP4})$ has a unique solution on $[-h,h]$.
\hfill$\blacksquare$
\end{dem}
The following corollary gives us the sufficient conditions for the existence and uniqueness of the solutions of the $q$-Cauchy problem $~(\ref{IP5})$.
\begin{corollary}\label{cor11}
Assume the functions $a_{j}(x):I\rightarrow \mathbb{C}$, $j=0,1,2$, and $b(x)\,:\,I \rightarrow \mathbb{X}$ satisfy the following conditions:
\begin{description}
  \item[(i)]  $a_{j}(x)$, $j=0,1,2$, and $b(x)$ are continuous at zero with $a_{0}(x)\neq0$ for all $x \in I$,
  \item[(ii)]  $\dfrac{a_{j}(x)}{a_{0}(x)}$ is bounded on $I$, $j=1,2$. Then
\begin{equation}\label{orderff'}
\left\{
  \begin{array}{cc}
      q\,a_0(x)\,\partial_q^2y(qx)\,+\,a_1(x)\,\partial_qy(x)\,+\,a_2(x)y(x)&\;=\;b(x),\quad \hbox{if $y$ is odd;}\\
      q\,a_0(x)\partial_q^2y(qx)\,+\,q\,a_1(x)\partial_qy(qx)\,+\,a_2(x)y(x)&\;=\;b(x),\quad \hbox{if $y$ is even,}
 \end{array}
\right.
    \end{equation}
with the initial conditions
\begin{equation}\label{abs2}
   \partial_{q}^{i-1}y(0)=b_i;\quad b_i\in\mathbb{C},\; i=1,2
 \end{equation}
\end{description}
has a unique solution on subinterval $J\subseteq I$ containing zero.
\end{corollary}
\begin{dem}
Dividing $~(\ref{orderff'})$ by $a_0(x)$, we get
\begin{equation}\label{IP6}
\left\{
  \begin{array}{cc}
      \partial_{q}^{2}y(qx)\;=\; A_{1}(x)\partial_{q}y(x)+A_{2}(x)y(x)+B(x),\quad \hbox{if $y$ is odd;}\\
      \partial_{q}^{2}y(qx)\;=\; q\,A_{1}(x)\partial_{q}y(qx)+A_{2}(x)y(x)+B(x),\quad \hbox{if $y$ is even,}
 \end{array}
\right.
    \end{equation}
where $A_{j}(x)=-q^{-1}\,\dfrac{a_{j}(x)}{a_{0}(x)}$ and $B(x)=\dfrac{b(x)}{a_{0}(x)}$. Since $A_{j}(x)$ and $B(x)$ are continuous at zero, the function $f(x,y_{1},y_{2})$, defined by
$$\left\{
  \begin{array}{cc}
      f(x,y_{1},y_{2})\;=\; A_{1}(x)\partial_{q}y(x)+A_{2}(x)y(x)+B(x),\quad \hbox{if $y$ is odd;}\\
      f(x,y_{1},y_{2})\;=\; q\,A_{1}(x)\partial_{q}y(qx)+A_{2}(x)y(x)+B(x),\quad \hbox{if $y$ is even,}
 \end{array}
\right.$$
is continuous at zero. Furthermore, $A_{j}(x)$ is bounded on $I$. Consequently, there is $L>0$ such that $|A_{j}(x)|\leq L$ for all $x \in I$. We can see that $f$ satisfies the Lipschitz condition with Lipschitz constant $L$. Thus, $f(x,y_{1},y_{2})$ satisfies the conditions of Corollary $~\ref{cor1}$. Hence, there exists a unique solution of $~(\ref{IP6})$ on $J$.
\hfill$\blacksquare$
\end{dem}
\section{Second Order Homogeneous Linear $q$-Difference Equations Associated with the Rubin's $q$-Difference Operator $\partial_q$}
Consider the second order non-homogeneous $q$-difference equation associated with the Rubin's $q$-difference operator $\partial_q$ in a neighborhood of zero
\begin{equation}
\left\{
  \begin{array}{cc}
      q\,a_0(x)\,\partial_q^2y(qx)\,+\,a_1(x)\,\partial_qy(x)\,+\,a_2(x)y(x)&\;=\;b(x),\quad \hbox{if $y$ is odd;}\\
      q\,a_0(x)\partial_q^2y(qx)\,+\,q\,a_1(x)\partial_qy(qx)\,+\,a_2(x)y(x)&\;=\;b(x),\quad \hbox{if $y$ is even,}
 \end{array}
\right.
    \end{equation}
with the initial conditions
\begin{equation}
   \partial_{q}^{i-1}y(0)= b_{i};\quad b_{i} \in{\mathbb{C}},\; i=1,2
 \end{equation}
where $a_i$, $i=0,1,2$, and $b$ are defined, continuous at zero and bounded on an interval $I$ containing zero such that $a_0(x)\neq 0$ for all $x\in I$. In this section, we shall study the second order homogeneous linear $q$-difference problem associated with the Rubin's $q$-difference operator $\partial_q$ of the form
\begin{equation}\label{order}
\left\{
  \begin{array}{cc}
      q\,a_0(x)\,\partial_q^2y(qx)\,+\,a_1(x)\,\partial_qy(x)\,+\,a_2(x)y(x)&\;=\;0,\quad \hbox{if $y$ is odd;}\\
      q\,a_0(x)\partial_q^2y(qx)\,+\,q\,a_1(x)\partial_qy(qx)\,+\,a_2(x)y(x)&\;=\;0,\quad \hbox{if $y$ is even.}
 \end{array}
\right.
    \end{equation}
The following result summarizes some properties of $~(\ref{order})$ which we can state at once.
\begin{prop}\label{pr}
Let $x$ in a subinterval $J$ of $I$ which contains zero. Then
\begin{enumerate}
  \item If $\phi_1$ and $\phi_2$ are two solutions of $~(\ref{order})$ and have the same parity, then
$$\phi(x)\;=\;c_1\phi_1(x)\,+\,c_2\phi_2(x),$$
where $c_1$ and $c_2$ are constants, is also a solution of $~(\ref{order})$.
  \item If $\phi$ is a solution of $~(\ref{order})$ such that $\partial_q^{\,i-1}\phi(0)\,=\,0$, $i=1,2$, then
$$\phi(x)\;=\;0.$$
\end{enumerate}
\end{prop}
\begin{dem}
\begin{enumerate}
  \item Since $\phi_1$ and $\phi_2$ are two solutions of Equation $~(\ref{order})$ and have the same parity, we distinguish two cases:
      \begin{description}
\item[--]   If $\phi_1$ and $\phi_2$ are odd, then, $\phi(x)\,=\,c_1\phi_1(x)\,+\,c_2\phi_2(x)$ is an odd function. Therefore, using the basic rules for $q$-differentiation, we get
  \begin{align*}
     q\,a_0(x)\,&\partial_q^{\,2}\phi(qx)\,+\,a_1(x)\,\partial_q\phi(x)\,+\,a_2(x)\,\phi(x)\\
\,=\;&q\,a_0(x)\partial_q^{\,2}\left[c_1\,\phi_1(qx)+c_2\,\phi_2(qx)\right]\,+\,a_1(x)\partial_q\left[c_1\,\phi_1(x)+c_2\,\phi_2(x)\right]\\
&\,+\;a_2(x)[c_1\,\phi_1(x)+c_2\,\phi_2(x)]\\
\,=\;&q\,a_0(x)\left[c_1\,\partial_q^{\,2}\phi_1(qx)+c_2\,\partial_q^2\phi_2(qx)\right]+a_1(x)\left[c_1\,\partial_q\phi_1(x)+c_2\,\partial_q\phi_2(x)\right]\\
&\,+\;a_2(x)[c_1\,\phi_1(x)+c_2\,\phi_2(x)]\\
\,=\;&c_1\,\left[q\,a_0(x)\partial_q^{\,2}\phi_1(qx)\,+\,a_1(x)\partial_q\phi_1(x)\,+\,a_2(x)\phi_1(x)\right]\\
&\,+\;c_2\,\left[q\,a_0(x)\partial_q^{\,2}\phi_2(qx)\,+\,a_1(x)\partial_q\phi_2(x)\,+\,a_2(x)\phi_2(x)\right]\\
\,=\;&0.
  \end{align*}
  \item[--]   If $\phi_1$ and $\phi_2$ are even, then, $\phi(x)\,=\,c_1\,\phi_1(x)\,+\,c_2\,\phi_2(x)$ is an even function. Therefore, using the basic rules for $q$-differentiation, we get
  \begin{align*}
     q\,a_0(x)\,&\partial_q^{\,2}\phi(qx)\,+\,q\,a_1(x)\,\partial_q\phi(qx)\,+\,a_2(x)\,\phi(x)\\
\,=\;&q\,a_0(x)\partial_q^{\,2}\left[c_1\,\phi_1(qx)+c_2\,\phi_2(qx)\right]\,+\,q\,a_1(x)\partial_q\left[c_1\,\phi_1(qx)+c_2\,\phi_2(qx)\right]\\
&\,+\;a_2(x)[c_1\,\phi_1(x)+c_2\,\phi_2(x)]\\
\,=\;&q\,a_0(x)\left[c_1\,\partial_q^{\,2}\phi_1(qx)+c_2\,\partial_q^2\phi_2(qx)\right]+q\,a_1(x)\left[c_1\,\partial_q\phi_1(qx)+c_2\,\partial_q\phi_2(qx)\right]\\
&\,+\;a_2(x)[c_1\,\phi_1(x)+c_2\,\phi_2(x)]\\
\,=\;&c_1\,\left[q\,a_0(x)\partial_q^{\,2}\phi_1(qx)\,+\,q\,a_1(x)\partial_q\phi_1(qx)\,+\,a_2(x)\phi_1(x)\right]\\
&\,+\;c_2\,\left[q\,a_0(x)\partial_q^{\,2}\phi_2(qx)\,+\,q\,a_1(x)\partial_q\phi_2(qx)\,+\,a_2(x)\phi_2(x)\right]\\
\,=\;&0.
\end{align*}
\end{description}
  Thus, $\phi(x)$ is also a solution of Equation $~(\ref{order})$.
  \item The function $\phi_0(x)$ which is identically zero in $I$ clearly satisfies $~(\ref{order})$ and the initial conditions
  $$\partial_q^{\,i-1}\phi_0(0)\;=\;0,\quad i=1,2.$$
  Thus $\phi(x)$ and $\phi_0(x)$ satisfy the same initial conditions at zero and therefore, by Corollary $~\ref{cor11}$, there exists a unique solution of $~(\ref{order})$. Then, we have
  $$\phi(x)\;=\;\phi_0(x)\;=\;0,$$
\end{enumerate}
for all $x$ in $I$.\hfill$\blacksquare$
\end{dem}
\begin{defn}
 A set of two linearly independent solutions of $~(\ref{order})$ is called a fundamental set of it.
 \end{defn}
The existence of fundamental sets of $~(\ref{order})$ is established in the following lemma.
\begin{lemma}\label{lem}
Let $b_{ij}$, $i,j=1,2$, be any real or complex numbers. For each $j=1,2$, let $\phi_j$ be the solution of $~(\ref{order})$ which satisfies the initial conditions
\begin{equation}
  \partial_q^{\,i-1}\phi_j(0)\;=\;b_{ij},\quad i=1,2.
\end{equation}
Then a necessary and sufficient condition that $\{\phi_1,\phi_2\}$ is a fundamental set of $~(\ref{order})$ is that $\det (b_{ij})\neq 0$.
\end{lemma}
\begin{dem}
\underline{Necessary:}  Let $\{\phi_1,\phi_2\}$ is a fundamental set but suppose that $\det (b_{ij})  = 0$. Then there are numbers $\alpha_j$, $j=1,2$, not all zero, such that
\begin{equation}\label{xx}
\left\{
  \begin{array}{ll}
    \alpha_1b_{11}\,+\,\alpha_2b_{12}&\;=\;0\\
    \alpha_1b_{21}\,+\,\alpha_2b_{22}&\;=\;0.
  \end{array}
\right.\end{equation}
Now define $\phi(x)\;=\;\alpha_1\phi_{1}(x)\,+\,\alpha_2\phi_{2}(x)$, for all $x$ in a subinterval $J$ of $I$ which contains zero.
By Proposition $~\ref{pr}$, $\phi$ is a solution of $~(\ref{order})$ and, by $~(\ref{xx})$, we have
$$\phi(0)\;=\;\partial_q\phi(0)\;=\;0.$$
Hence, by Proposition $~\ref{pr}$, $\phi(x)\,=\,0$. But, since the $\alpha_j$, $j=1,2$, are not all zero, this contradicts the linear independence of the $\phi_j(x)$, and so we must have $\det (b_{ij})\neq 0$.\\
\underline{Sufficient:}  Let $\det (b_{ij})\neq 0$. Then we have to show that the relation
\begin{equation}\label{rr}\alpha_1\phi_{1}(x)\,+\,\alpha_2\phi_{2}(x)\;=\;0\end{equation}
is possible for all $x$ in $J$ only when $\alpha_1\,=\,\alpha_2\,=\,0$. Differentiating $(i-1)$ times $~(\ref{rr})$, $i=1,2$, and putting $x=0$, we obtain the equations $~(\ref{xx})$. Since $\det (b_{ij})\neq 0$, $~(\ref{xx})$ implies that $\alpha_1\,=\,\alpha_2\,=\,0$, as required.
\hfill$\blacksquare$ \end{dem}
 \section{The $q$-Wronskian Associated with the Rubin's $q$-Difference Operator $\partial_q$}
 To determine if two solutions of the second order $q$-difference equations associated with the Rubin's $q$-difference operator $\partial_q$ form a fundamental set, we introduce a $q$-analogue of the Wronskian. 
\begin{defn}\label{wq}
  Let  $y_1$ and $y_2$ be two $q$-differentiable functions defined on a $q$-geometric set $\mathcal{A}$. The $q$-Wronskian associated with the Rubin's $q$-difference operator $\partial_q$ of the functions $y_1$ and $y_2$ which will be denoted by $W_q(y_1,y_2)$ is defined by:
  \begin{description}
  \item[--]  If $y_1$ is even and $y_2$ is odd, we have
  \begin{equation}\label{ff1}
  W_q(y_1,y_2)(x)\;=\;{\displaystyle \begin{vmatrix}y_{1}(x)&y_{2}(x)\\q\,\partial_qy_{1}(qx)&\partial_qy_{2}(x) \end{vmatrix}}\;=\;y_{1}(x)\partial_qy_{2}(x)\,-q\,\,y_{2}(x)\partial_qy_{1}(qx).
\end{equation}
 \item[--]  If $y_1$ and $y_2$ are odd, we have
  \begin{equation}\label{ff2}
  W_q(y_1,y_2)(x)\;=\;{\displaystyle \begin{vmatrix}y_{1}(x)&y_{2}(x)\\\partial_qy_{1}(x)&\partial_qy_{2}(x) \end{vmatrix}}\;=\;y_{1}(x)\partial_qy_{2}(x)\,-\,y_{2}(x)\partial_qy_{1}(x).
\end{equation}
 \item[--]  If $y_1$ and $y_2$ are even, we have
  \begin{equation}\label{ff3}
  W_q(y_1,y_2)(x)\;=\;{\displaystyle \begin{vmatrix}y_{1}(x)&y_{2}(x)\\q\,\partial_qy_{1}(qx)&q\,\partial_qy_{2}(qx) \end{vmatrix}}
  \;=\;q\,y_{1}(x)\partial_qy_{2}(qx)\,-\,q\,y_{2}(x)\partial_qy_{1}(qx).
\end{equation}
\end{description}
It is easy to see that, if $q$ tends to $1^{-}$, the $q$-Wronskian tends to the ordinary Wronskian
\begin{equation}\label{w}
 W(y_1,y_2)(x)\;=\;{\displaystyle \begin{vmatrix}y_{1}(x)&y_{2}(x)\\y_{1}\,'(x)&y_{2}\,'(x) \end{vmatrix}}\;=\;y_{1}(x)y_{2}'(x)\,-\,y_{2}(x)y_{1}'(x).
\end{equation}
\end{defn}
\begin{corollary}
  Let $y_1$ and $y_2$ be two differentiable functions defined on a $q$-geometric set $\mathcal{A}$. 
  The $q$-Wronskian which is defined by $~(\ref{wq})$ can be rewritten in the form
\begin{equation}\label{n2}
  W_q(y_1,y_2)(x)\;=\;\dfrac{y_{2}(x)y_{1}(qx)\,-\,y_{1}(x)y_{2}(qx)}{(1-q)x},\quad x\in\mathcal{A}\backslash\{0\}.
\end{equation}
\end{corollary}
\begin{prop}\label{pl}
Let $y_1, y_2$ be two solutions of $~(\ref{order})$ defined on $\mathcal{A}$. Then, for all $x\in\mathcal{A}\backslash\{0\}$, we have
\begin{description}
  \item[--]  If $y_1$ and $y_2$ have opposite parity, then
  \begin{equation}\label{pl1}
 \partial_qW_q(y_1,y_2)(x)\;=\;y_1(x)\,\partial_q^2\,y_2(x)\,-\,y_2(x)\,\partial_q^2\,y_1(x).
\end{equation}
  \item[--]  If $y_1$ and $y_2$ have the same parity, then
  \begin{equation}\label{pl1'}
  \partial_qW_q(y_1,y_2)(x)\;=\;q\,y_1(qx)\,\partial_q^2\,y_{2}(qx)\,-\,q\,y_2(qx)\,\partial_q^2\,y_{1}(qx).
\end{equation}
\end{description}
\end{prop}
\begin{dem}
To prove $~(\ref{pl1})$, we involve several cases such as:
\begin{description}
\item[--]  If $y_1$ is an even function and $y_2$ is an odd function, the $q$-Wronskian is given by $~(\ref{ff1})$. Using the fact that $\partial_qy_1$ is odd and $\partial_qy_2$ is even, we get by the help of $~(\ref{marwa2})$ and $~(\ref{marwa3})$ as follows
\begin{align*}
\partial_qW_q(y_1,y_2)(x)&\;=\;\partial_q\left[y_{1}(x)\partial_qy_{2}(x)\right]-q\,\partial_q\left[y_{2}(x)\partial_qy_{1}(qx)\right]\\
&\;=\;\partial_q\,y_2(q^{-1}x)\,\partial_qy_1(x)\,+\,y_1(x)\,\partial_q^2y_2(x)\\
&\qquad\;-\,q\,\left[q^{-1}\,\partial_q\,y_1(x)\,\partial_qy_2(q^{-1}x)\,+\,q^{-1}\,y_2(x)\,\partial_q^2y_1(x)\right]\\
 &\;=\;y_1(x)\,\partial_q^2y_2(x)\,-\,y_2(x)\,\partial_q^2y_1(x).
 \end{align*}
\item[--]  If $y_1$ and $y_2$ are odd, the $q$-Wronskian is given by $~(\ref{ff2})$. Using the fact that $\partial_qy_1$ and $\partial_qy_2$ are even, we get by the help of $~(\ref{marwa1})$ as follows
 \begin{align*}
 \partial_qW_q(y_1,y_2)(x)&\;=\;\partial_q\left[y_{1}(x)\partial_qy_{2}(x)\right]\,-\,\partial_q\left[y_{2}(x)\partial_qy_{1}(x)\right]\\
&\;=\;q\,y_1(qx)\,\partial_q^2y_{2}(qx)\,-\,q\,y_2(qx)\,\partial_q^2y_{1}(qx).
\end{align*}
\item[--]  If $y_1$ and $y_2$ are even, the $q$-Wronskian is given by $~(\ref{ff3})$. Using the fact that $\partial_qy_1$ and $\partial_qy_2$ are odd, we get by the help of $~(\ref{marwa1})$ as follows
 \begin{align*}
 \partial_qW_q(y_1,y_2)(x)&\;=\;q\,\partial_q\left[y_{1}(x)\partial_qy_{2}(qx)\right]\,-\,q\,\partial_q\left[y_{2}(x)\partial_qy_{1}(qx)\right]\\
&\;=\;q\,y_1(qx)\,\partial_q^2y_{2}(qx)\,-\,q\,y_2(qx)\,\partial_q^2y_{1}(qx).
\end{align*}
\end{description}
This completes the proof.
\hfill$\blacksquare$ \end{dem}
\begin{theorem}{\bf($q$-analogue of Abel's theorem)}\label{Abel}
  Let $I$ be an interval containing zero. If $y_1, y_2$ are solutions of $~(\ref{order})$ in a subinterval $J$ of $I$, $J=[-h,h]$, $h > 0$, then their $q$-Wronskian associated with the Rubin's $q$-difference operator $\partial_q$ satisfies the linear first order $q$-difference equation as follows:
  \begin{description}
  \item[--]  If $y_1$ and $y_2$ have opposite parity, then
  \begin{equation}\label{Wqs1}
    \partial_q\,W_q(y_1,y_2)(x)\;=\;-q^{-1}\,E(q^{-1}x)\,W_q(y_1,y_2)(q^{-1}x);
  \end{equation}
  \item[--]  If $y_1$ and $y_2$ have the same parity, then
 \begin{equation}\label{Wqs2}
    \partial_q\,W_q(y_1,y_2)(x)\;=\;-E(x)\,W_q(y_1,y_2)(x);
  \end{equation}
\end{description}
  where
  \begin{equation}\label{R(x)}
    E(x)\;=\;\dfrac{a_1(x)\,+\,x(1-q)\,{a_2(x)}}{a_0(x)},
  \end{equation}
  for all $x\in J\backslash\{0\}$.
\end{theorem}
\begin{dem}
To prove $~(\ref{Wqs1})$ and $~(\ref{Wqs2})$, we involve several cases according to the parity of $y_1$ and $y_2$. Then, we have
\begin{description}
\item[--]  If $y_1$ is an even function and $y_2$ is an odd function, from $~(\ref{order})$, we obtain
$$\partial_q^2y_{1}(x)\;=\;-\,\dfrac{q\,a_1(q^{-1}x)\partial_qy_1(x)+a_2(q^{-1}x)y_1(q^{-1}x)}{q\,a_0(q^{-1}x)}$$
and
$$\partial_q^2y_{2}(x)\;=\;-\dfrac{a_1(q^{-1}x)\partial_qy_2(q^{-1}x)+a_2(q^{-1}x)y_2(q^{-1}x)}{q\,a_0(q^{-1}x)}.$$
It follows from $~(\ref{pl1})$ that
\begin{align*}
 \partial_qW_q(y_1,y_2)(x)\;=&\;y_1(x)\,\partial_q^2y_{2}(x)\,-\,y_2(x)\,\partial_q^2y_{1}(x)\\
 \;=&\;-q^{-1}\,\dfrac{a_1(q^{-1}x)}{a_0(q^{-1}x)}\left(y_1(x)\,\partial_qy_2(q^{-1}x)\,-\,q\,y_2(x)\,\partial_qy_1(x)\right)\\
  &\;-\;q^{-1}\,\dfrac{a_2(q^{-1}x)}{a_0(q^{-1}x)}\left(y_1(x)\,y_2(q^{-1}x)\,-\,y_2(x)\,y_1(q^{-1}x)\right).
  \end{align*}
Hence, using $~(\ref{ff1})$ and $~(\ref{n2})$, we obtain
 $$\partial_qW_q(y_1,y_2)(x)\;=\;-\,q^{-1}\,\left[\dfrac{a_1(q^{-1}x)\,+\,q^{-1}x(1-q)\,{a_2(q^{-1}x)}}{a_0(q^{-1}x)}\right]\,W_q(y_1,y_2)(q^{-1}x).$$
From $~(\ref{R(x)})$, we conclude that
    $$\partial_qW_q(y_1,y_2)(x)\;=\;-\,q^{-1}\,E(q^{-1}x)\,W_q(y_1,y_2)(q^{-1}x),\quad x\in J\backslash\{0\}.$$
 \item[--]  If $y_1$ and $y_2$ are odd, from $~(\ref{order})$, we obtain
 $$\partial_q^2y_{1}(qx)\;=\;-\dfrac{a_1(x)\partial_qy_1(x)+a_2(x)y_1(x)}{q\,a_0(x)}$$
 and
 $$\partial_q^2y_{2}(qx)\;=\;-\dfrac{a_1(x)\partial_qy_2(x)+a_2(x)y_2(x)}{q\,a_0(x)}.$$
 It follows from $~(\ref{pl1'})$ that
\begin{align*}
 \partial_qW_q(y_1,y_2)(x)\;=&\;q\,y_1(qx)\,\partial_q^2y_{2}(qx)\,-\,q\,y_2(qx)\,\partial_q^2y_{1}(qx)\\
 \;=&\;-\dfrac{a_1(x)}{a_0(x)}\left(y_1(qx)\,\partial_qy_2(x)\,-\,y_2(qx)\,\partial_qy_1(x)\right)\\
  &\;-\;\dfrac{a_2(x)}{a_0(x)}\left(y_1(qx)\,y_2(x)\,-\,y_2(qx)\,y_1(x)\right).
  \end{align*}
Hence, using $~(\ref{ff2})$ and $~(\ref{n2})$, we obtain
 $$\partial_qW_q(y_1,y_2)(x)\;=\;-\,\left[\dfrac{a_1(x)\,+\,(1-q)\,x\,{a_2(x)}}{a_0(x)}\right]\,W_q(y_1,y_2)(x).$$
From $~(\ref{R(x)})$, we conclude that
    $$\partial_qW_q(y_1,y_2)(x)\;=\;-\,E(x)\,W_q(y_1,y_2)(x),\quad x\in J\backslash\{0\}.$$
\item[--]  If $y_1$ and $y_2$ are even, from $~(\ref{order})$, we obtain
 $$\partial_q^2y_{1}(qx)\;=\;-\dfrac{q\,a_1(x)\partial_qy_1(qx)+a_2(x)y_1(x)}{q\,a_0(x)}$$
 and
 $$\partial_q^2y_{2}(qx)\;=\;-\dfrac{q\,a_1(x)\partial_qy_2(qx)+a_2(x)y_2(x)}{q\,a_0(x)}.$$
 It follows from $~(\ref{pl1'})$ that
\begin{align*}
 \partial_qW_q(y_1,y_2)(x)\;=&\;q\,y_1(qx)\,\partial_q^2y_{2}(qx)\,-\,q\,y_2(qx)\,\partial_q^2y_{1}(qx)\\
 \;=&\;-\dfrac{a_1(x)}{a_0(x)}\left(q\,y_1(qx)\,\partial_qy_2(qx)\,-\,q\,y_2(qx)\,\partial_qy_1(qx)\right)\\
  &\;-\;\dfrac{a_2(x)}{a_0(x)}\left(y_1(qx)\,y_2(x)\,-\,y_2(qx)\,y_1(x)\right).
  \end{align*}
Hence, using $~(\ref{ff3})$ and $~(\ref{n2})$, we obtain
 $$\partial_qW_q(y_1,y_2)(x)\;=\;-\,\left[\dfrac{a_1(x)\,+\,(1-q)\,x\,{a_2(x)}}{a_0(x)}\right]\,W_q(y_1,y_2)(x).$$
From $~(\ref{R(x)})$, we conclude that
    $$\partial_qW_q(y_1,y_2)(x)\;=\;-\,E(x)\,W_q(y_1,y_2)(x),\quad x\in J\backslash\{0\}.$$
\end{description}
This finishes the proof of the theorem.
\hfill$\blacksquare$ \end{dem}
The following theorem gives a $q$-type Liouville’s formula for the $q$-Wronskian associated with the Rubin’s $q$-difference operator $\partial_q$.
\begin{theorem}{\bf($q$-Liouville’s formula)}\label{li}
Suppose that $x(1\,-\,q)\,E(x)\neq -1$ for all $x$ in a subinterval $J$ of $I$ which contains zero. Then the $q$-Wronskian of any set of solutions $\{y_1,y_2\}$ of $~(\ref{order})$ is given by
\begin{equation}\label{Liou}
W_q(y_1,y_2)(x)\;=\;\dfrac{W_q(y_1,y_2)(0)}{\D\prod_{k=0}^{+\infty}\left(1\,+\,x(1\,-\,q)q^k\,E(xq^k)\right)},\quad x\in J
\end{equation}
where $E(x)$ is defined by $~(\ref{R(x)})$.
\end{theorem}
\begin{dem}
Suppose that $x(1\,-\,q)\,E(x)\neq -1$ for all $x$ in a subinterval $J$ of $I$ which contains zero. Here we distinguish two cases:
\begin{description}
\item[--] In case $y_1$ and $y_2$ have opposite parity, the $q$-Wronskian $~(\ref{ff1})$ is even. Using Theorem$~\ref{Abel}$, the $q$-Wronskian $W_q(y_1,y_2)(x)$ satisfies the linear first order $q$-difference equation $~(\ref{Wqs1})$. Then, we get
$$W_q(y_1,y_2)(x)\;=\;(1\,+\,q^{-1}\,x(1\,-\,q)\,E(q^{-1}x))\,W_q(y_1,y_2)(q^{-1}x),\quad x\in J\backslash\{0\}.$$
Replacing $x$ by $qx$ in the previous equation, we obtain
$$W_q(y_1,y_2)(qx)\;=\;(1\,+\,x(1\,-\,q)\,E(x))\,W_q(y_1,y_2)(x),\quad x\in J\backslash\{0\}.$$
\item[--] In case $y_1$ and $y_2$ have the same parity, the $q$-Wronskian is odd. Using Theorem $~\ref{Abel}$, the $q$-Wronskian $W_q(y_1,y_2)(x)$ satisfies the linear first order $q$-difference equation $~(\ref{Wqs2})$. Then, we get
 $$W_q(y_1,y_2)(qx)\;=\;(1\,+\,x(1\,-\,q)\,E(x))\,W_q(y_1,y_2)(x),\quad x\in J\backslash\{0\},$$
\end{description}
Hence, under the assumption $1\,+\,x(1\,-\,q)\,E(x)\,\neq\, 0$, we deduce for all solutions $y_1$ and $y_2$ that
$$W_q(y_1,y_2)(x)\;=\;\dfrac{W_q(y_1,y_2)(qx)}{1\,+\,x(1\,-\,q)\,E(x)}.$$
Now we use induction on $n$ to see that we have for $n\in\mathbb{N} $,
$$W_q(y_1,y_2)(x)\;=\;\dfrac{W_q(y_1,y_2)(xq^n)}{\D\prod_{k=0}^{n-1}\left(1\,+\,x(1\,-\,q)q^k\,E(xq^k)\right)}.$$
Since all functions $\dfrac{a_j}{a_0}$, $1\leq j\leq n$, are continuous at zero, then $\D\sum_{k=0}^{+\infty}q^k\,|E(xq^k)|$ is convergent. Consequently,
$$\D\prod_{k=0}^{n-1}\left(1\,+\,x(1\,-\,q)q^k\,E(xq^k)\right)\quad \text{converges for every $x\in J$}.$$
Letting $n \rightarrow \infty$ and noting that $0 < q < 1$, one gets $q^n$ tends to 0. Thus, using the continuity of $W_q(y_1,y_2)(x)$ at zero, $~(\ref{Liou})$ follows.
\hfill$\blacksquare$ \end{dem}
\begin{corollary}\label{bb5}
Let $\{y_1,y_2\}$ be a set of solutions of $~(\ref{order})$ in some subinterval $J$ of $I$ which contains zero. Then $W_q(y_1,y_2)(x)$ is either never zero in $I$ if and only if $\{y_1,y_2\}$ is a fundamental set of $~(\ref{order})$.
\begin{dem}
From Lemma $~\ref{lem}$, the functions $\{y_1,y_2\}$ form a fundamental set of $~(\ref{order})$ if and only if $W_q(0)\,\neq\,0$. Hence, the result is a direct consequence of Theorem $~\ref{li}$.
\hfill$\blacksquare$ \end{dem}
\end{corollary}
\section{Homogeneous Equations with constant coefficients}
We start now to study the second order homogeneous linear $q$-difference equations associated with the Rubin's $q$-difference operator that contain constant coefficients only:
\begin{equation}\label{order"}
a\,\partial_{q}^{2}y(x)\,+\,b\,y(x)=\,0,\quad a \neq 0,\end{equation}
with the initial conditions
\begin{equation}\label{vb}
y_1(0)\,=\,1,\quad \partial_qy_1(0)\,=\,0\quad \text{and} \quad y_2(0)\,=\,0,\quad \partial_qy_2(0)\,=\,1,\end{equation}
respectively, where $a$ and $b$ are constants. In the following result, we will solve the equation $~(\ref{order"})$.
\begin{prop}\label{fgf'}
 Let $x\in\mathbb{R}_q$. Then, the $q$-difference equation $~(\ref{order"})$ have two distinct solutions
$$y_1(x)\;=\;\cos(\sqrt{\dfrac{b}{a}}x,q^2)\quad\text{and}\quad y_2(x)\;=\;\sqrt{\dfrac{a}{b}}\sin(\sqrt{\dfrac{b}{a}}x,q^2),$$
respectively, where $a$ and $b$ are constants, $a \neq 0$.
\end{prop}
\begin{dem}
First, it easy to see that, if the solution $y(x)$ of the $q$-difference problem $~(\ref{order"})$ is an analytic function at the origin, it can be developed in entire powers series. Then we have
\begin{equation}
  y(x)\;=\;\D\sum_{n=0}^{\infty}a_n\,x^n,\quad x\in\mathbb{R}_q.
\end{equation}
Here we distinguish two cases:
\begin{description}
  \item[--] If $n\,=\,2p$, $p\in\mathbb{N}$, the solution $y_1$ of $~(\ref{order"})$ is even. Then, we have
\begin{equation}\label{fg1}
 \partial_q^2y_1(x)\;=\;\D\sum_{p\geq0}^{}\,q^{-2(p+1)}\dfrac{1-q^{2p+2}}{1-q}\,\dfrac{1-q^{2p+1}}{1-q}\,a_{2p+2}\,x^{2p},\quad p\in\mathbb{N}.
\end{equation}
Loading $~(\ref{fg1})$ in $~(\ref{order"})$ and equating the coefficients of $x^{2p}$, one obtains
$$a\,q^{-2(p+1)}\,\dfrac{1-q^{2p+2}}{1-q}\,\dfrac{1-q^{2p+1}}{1-q}\,a_{2p+2}\;=\;-b\,a_{2p},\quad p\in\mathbb{N}.$$
Hence, We get the following recurrence relation for the coefficients:
$$a_{2p}\;=\;-\,\dfrac{b}{a}\,q^{2p}\,\dfrac{(1-q)^2}{(1-q^{2p})\,(1-q^{2p-1})}\,a_{2p-2},\quad p\geq 1.$$
So, by induction on $p$ and the fact that $ a_{0}\,=\,1$, we obtain:
$$a_{2p}\;=\;(-1)^p\,q^{p(p+1)}\,\left(\dfrac{b}{a}\right)^p\,\dfrac{(1-q)^{2p}}{(q,q)_{2p}}\;=\;(-1)^p\,\dfrac{q^{p(p+1)}}{[2p]_q!}\,\,\left(\sqrt{\dfrac{b}{a}}\right)^{2p}\quad p\in\mathbb{N}.$$
From the definitions $~(\ref{cos})$, we get
$$y_1(x)\;=\;\cos(\sqrt{\dfrac{b}{a}}x,q^2).$$
\item[--] If $n\,=\,2p+1$, $p\in\mathbb{N}$, the solution $y_2$ of $~(\ref{order"})$ is odd. Then, we have
\begin{equation}\label{fg4}
  \partial_{q}^{2}y_2(x)\;=\;\D\sum_{p\geq0}^{}\,q^{-2(p+1)}\dfrac{1-q^{2p+3}}{1-q}\,\dfrac{1-q^{2p+2}}{1-q}\,a_{2p+3}\,x^{2p+1},\quad p\in\mathbb{N}.
\end{equation}
Loading $~(\ref{fg4})$ in $~(\ref{order"})$ and equating the coefficients of $x^{2n+1}$, one obtains
$$a\,q^{-2(p+1)}\,\dfrac{1-q^{2p+3}}{1-q}\,\dfrac{1-q^{2p+2}}{1-q}\,a_{2p+3}\;=\;-b\,a_{2p+1},\quad p\in\mathbb{N}.$$
Hence, we get the following recurrence relation for the coefficients:
$$a_{2p+1}\;=\;-\,\dfrac{b}{a}\,q^{2p}\,\dfrac{(1-q)^2}{(1-q^{2p+1})\,(1-q^{2p})}\,a_{2p-1},\quad p\geq 1.$$
So, by induction on $p$ and the fact that $ a_{1}\,=\,1$, we obtain:
\begin{align*}
a_{2p+1}\;=&\;(-1)^p\,q^{p(p+1)}\,\left(\dfrac{b}{a}\right)^p\,\dfrac{(1-q)^{2p+1}}{(q,q)_{2p+1}}\\
\;=&\;\sqrt{\dfrac{a}{b}}\,(-1)^p\,\dfrac{q^{p(p+1)}}{[2p+1]_q!}\,\left(\sqrt{\dfrac{b}{a}}\right)^{2p+1},\quad p\in\mathbb{N}.\end{align*}
From the definitions $~(\ref{sin})$, we get
$$y_2(x)\;=\;\sqrt{\dfrac{a}{b}}\sin(\sqrt{\dfrac{b}{a}}x,q^2).$$
\end{description}
Then, the functions $y_1(x)$ and $y_2(x)$, $x\in\mathbb{R}_q$, are solutions of $~(\ref{order"})$. This completes the proof.
\hfill$\blacksquare$
\end{dem}
\begin{corollary}\label{nb}
  Let $x\in\mathbb{R}_q$. Then Equation $~(\ref{order"})$ can be written another way:
  \begin{equation}\label{eex1'}
\left\{
  \begin{array}{cc}
      a\,\partial_q^2y(qx)\,-\,b(1-q)\,x\,\partial_qy(qx)\,+\,b\,y(x)\;=\;0,\quad \hbox{if $y$ is odd;}\\
      a\,\partial_q^2y(qx)\,-\,b\,q\,(1-q)\,x\,\partial_qy(qx)\,+\,b\,y(x)\;=\;0,\quad \hbox{if $y$ is even,}
 \end{array}
\right.
    \end{equation}
where $a$ and $b$ are constants, $a \neq 0$.
\end{corollary}
\begin{dem}
Two cases arise according to $y$ is even or odd. Then we have
\begin{description}
\item[--]  If $y$ is odd, the $\partial_qy$ would be an even function. Detailing the Rubin’s $q$-difference operator $\partial_q$ in $~(\ref{order"})$, the equation becomes
\begin{equation}\label{eex41}
a\,\partial_qy(q^{-1}x)\,-\,a\,\partial_qy(x)\,+\,b\,(1-q)\,x\,y(x)\;=\;0.\end{equation}
Applying the $\partial_q\,\Lambda_q$ derivative in the previous equation and using $~(\ref{marwa3})$, one gets
$$a\,\partial_q^2y(x)\,-\,a\,q\,\partial_q^2y(qx)\,+\,b\,(1-q)\,y(x)\,+\,b\,q\,(1-q)\,x\,\partial_qy(x)\;=\;0.$$
Loading $~(\ref{order"})$ in the last equation and next dividing the lefthand by $(-q)$, we obtain
$$a\,\partial_q^2y(qx)\,-\,b\,(1-q)\,x\,\partial_qy(x)\,+\,b\,y(x)\;=\;0.$$
  \item[--]  If $y$ is even, the $\partial_qy$ would be an odd function. Detailing the Rubin’s $q$-difference operator $\partial_q$ in $~(\ref{order"})$, the equation reads
\begin{equation}\label{eex4}
a\,\partial_qy(x)\,-\,a\,\partial_qy(qx)\,+\,b\,(1-q)\,x\,y(x)\;=\;0.\end{equation}
Applying the $\partial_q$ derivative in the previous equation and using $~(\ref{marwa1})$, one gets
$$a\,\partial_q^2y(x)\,-\,a\,q\,\partial_q^2y(qx)\,+\,b\,(1-q)y(x)\,+\,b\,q^2\,(1-q)\,x\,\partial_qy(qx)\;=\;0.$$
Loading $~(\ref{order"})$ in the last equation and next dividing the lefthand by $(-q)$, we obtain
$$a\,\partial_q^2y(qx)\,-\,b\,q\,(1-q)\,x\,\partial_qy(qx)\,+\,b\,y(x)\;=\;0.$$
\end{description}
Thus the result hold.
\hfill$\blacksquare$
\end{dem}
Let us now calculate the $q$-Wronskian of the solutions of $~(\ref{order"})$ in the following result:
\begin{prop}\label{nb1}
Let $x\in \mathbb{R}_q$. Then the $q$-Wronskian $W_q(x)$ of the solutions of the $q$-difference equation $~(\ref{order"})$ subject to the initial conditions $~(\ref{vb})$ is given by
\begin{equation}
  W_q(x)\;=\;1,\quad  x\in\mathbb{R}_q.
\end{equation}
Moreover, the set $\{y_1, y_2\}$ form a fundamental set of $~(\ref{order"})$.
\end{prop}
\begin{dem}
Comparing $~(\ref{eex1'})$ with $~(\ref{order})$, we get
$$a_0(x)\,=\,aq^{-1},\quad a_1(x)\,=\,-b\,x\,(1-q)\quad \text{and}\quad a_2(x)\,=\,b,$$
where $a$, $b$ are constants and $a\neq0$. Using the formula $~(\ref{R(x)})$, we obtain
$$E(x)\;=\;0,\quad x\in\mathbb{R}_q.$$
So, from the $q$-Liouville’s formula $~(\ref{Liou})$, we deduce that
$$W_q(x)\;=\;W_q(0),\quad x\in\mathbb{R}_q.$$
But, it follows from Proposition $~\ref{fgf'}$ that the solutions of the $q$-difference equation $~(\ref{order"})$ are the functions
$$y_1(x)\;=\;\cos(\sqrt{\dfrac{b}{a}}x,q^2)\quad\text{and}\quad y_2(x)\;=\;\sqrt{\dfrac{a}{b}}\sin(\sqrt{\dfrac{b}{a}}x,q^2),$$
with the initial conditions $~(\ref{vb})$, where $a$ and $b$ are constants, $a \neq 0$.
Then, by $~(\ref{ff1})$, we get
\begin{align*}
W_q(0)&\;=\;W_q(\cos(\sqrt{\dfrac{b}{a}}x,q^2),\sqrt{\dfrac{a}{b}}\sin(\sqrt{\dfrac{b}{a}}x,q^2))\mid_{x=0}\\
&\;=\;\left(\cos(\sqrt{\dfrac{b}{a}}x,q^2)\cos(\sqrt{\dfrac{b}{a}}x,q^2)
\,+\,q\,\sin(\sqrt{\dfrac{b}{a}}x,q^2)\sin(q\,\sqrt{\dfrac{b}{a}}x,q^2)\right)\mid_{x=0}\\
&\;=\;1.
\end{align*}
Consequently, $W_q(x)\,=\,1$, for all $x\in \mathbb{R}_q$. Moreover, by $~(\ref{bb5})$, the set $\{y_1, y_2\}$ form a fundamental set of $~(\ref{order"})$ and the result is proved.
\hfill$\blacksquare$
\end{dem}
\section{Examples}
\begin{ex} We purpose to solve the following $q$-difference equation associated with the Rubin’s $q$-difference operator
\begin{equation}\label{fg}
    \partial_q^2y(x)\,+\,y(x)\,=\,0,\quad x\in\mathbb{R}_q
\end{equation}
with the initial conditions
$$y_1(0)\,=\,1,\quad \partial_qy_1(0)\,=\,0\quad \text{and} \quad y_2(0)\,=\,0,\quad \partial_qy_2(0)\,=\,1,$$
respectively.\\
By Proposition$~\ref{fgf'}$, the solutions of $~(\ref{fg})$ are the functions $\cos(x,q^2)$, $\sin(x,q^2)$, $x\in\mathbb{R}_q$, respectively.\\
%
From Proposition $~\ref{nb1}$, the $q$-Wronskian $W_q(x)$ of the solutions of the $q$-difference equation $~(\ref{fg})$ is given by
$$W_q(x)\;=\;1,\quad x\in\mathbb{R}_q.$$
consequently, the set $\{\cos(x,q^2)$, $\sin(x,q^2)\}$ form a fundamental set of $~(\ref{fg})$.\end{ex}
 \begin{ex} We define a pair of basic $q$-trigonometric functions by studying the solutions of the
second order q-difference equation associated with the Rubin’s $q$-difference operator
\begin{equation}\label{fg0}
    q\,\partial_q^2y(x)\,+\,y(x)\,=\,0,\quad x\in\mathbb{R}_q
\end{equation}
with the initial conditions
$$y_1(0)\,=\,1,\quad \partial_qy_1(0)\,=\,0\quad \text{and} \quad y_2(0)\,=\,0,\quad \partial_qy_2(0)\,=\,1,$$
respectively.\\
By Proposition$~\ref{fgf'}$, the functions $\cos(q^{-1/2}x,q^2)$, $q^{1/2}\,\sin(q^{-1/2}x,q^2)$, $x\in\mathbb{R}_q$, are solutions of $~(\ref{fg0})$, respectively.\\
%
From Proposition $~\ref{nb1}$, the $q$-Wronskian $W_q(x)$ of the solutions of the $q$-difference equation $~(\ref{fg0})$ is given by
$$W_q(x)\;=\;1,\quad x\in\mathbb{R}_q.$$
consequently, the set $\{\cos(q^{-1/2}x,q^2)$, $q^{1/2}\,\sin(q^{-1/2}x,q^2)$ form a fundamental set of $~(\ref{fg0})$.\end{ex}

\end{document}